\documentclass[12pt,letterpaper]{article}
\usepackage{graphicx}
\usepackage{amsmath}
\usepackage{amsfonts}
\usepackage{amsthm}
\usepackage{amssymb}
\usepackage{color}

\usepackage{verbatim} 
\usepackage{hyperref}
\usepackage{breakurl}

\theoremstyle{remark}

\newcommand{\be}{\begin{equation}}
\newcommand{\ee}{\end{equation}}
\newcommand{\bea}{\begin{eqnarray}}
\newcommand{\eea}{\end{eqnarray}}
\newcommand{\beas}{\begin{eqnarray*}}
\newcommand{\eeas}{\end{eqnarray*}}

\begin{document}

\title{Martin Gardner's Mistake}

\author{Tanya Khovanova\\MIT}

\maketitle

\begin{abstract}
When Martin Gardner first presented the Two-Children problem, he made a mistake in its solution. Later he corrected the mistake in another publication, but unfortunately his incorrect solution is more widely known than his correction. In fact, a Tuesday-Child variation of this problem went viral in 2010, and the same flaw keeps reappearing in solutions for that problem as well. In this article, I would like to popularize Martin Gardner's correction and conduct a detailed discussion of the new problem.
\end{abstract}

\section{Martin Gardner}

Martin Gardner's writing is amazingly accurate and reliable. The fact that he made a mistake is simply a testimonial to the difficulty of the problem.

\subsection{Two-Children Problem}

The following problem appeared in Martin Gardner's \emph{Scientific American} column in 1959. It was later republished in his book \emph{The Second Scientific American Book of Mathematical Puzzles and Diversions} \cite{Gardner} in the chapter ``Nine More Problems:"

\begin{quote}
Mr. Smith has two children. At least one of them is a boy. What is the probability that both children are boys?

Mr. Jones has two children. The older child is a girl. What is the probability that both children are girls?
\end{quote}

There was never any controversy about Mr. Jones who with probability 1/2 has two girls, so we will just ignore him and his two kids. Instead, we will concentrate on Mr. Smith. Here's the solution that Martin Gardner published together with the problem:

\begin{quote}
If Smith has two children, at least on of which is a boy, we have three equally probable cases: Boy-boy, Boy-girl, Girl-boy. In only one case are both children boys, so the probability that both are boys is 1/3.
\end{quote}

\subsection{The Corrected Solution}

Later he wrote a column titled ``Probability and Ambiguity," which was republished in the same book as the original puzzle \cite{Gardner}. In this column he comes back to Mr. Smith, correcting himself:

\begin{quote}
... the answer depends on the procedure by which the information ``at least one is a boy" is obtained.
\end{quote}

He suggested two potential procedures.

\begin{enumerate}
\item Pick all the families with two children, one of which is a boy. If Mr. Smith is chosen randomly from this list, then the answer is 1/3. 
\item Pick a random family with two children; suppose the father is Mr. Smith. Then if the family has two boys, Mr. Smith says, ``At least one of them is a boy." If he has two girls, he says, ``At least one of them is a girl." If he has a boy and a girl he flips a coin to say one or another of those two sentences. In this case the probability that both children are the same sex is 1/2.
\end{enumerate}

Thus, the original problem without a specified procedure is ambiguous.

\subsection{More Procedures}

I call the first procedure ``boy-centered," because from the start we know that we are talking about boys. Correspondingly, the second procedure is called ``gender neutral."

Martin Gardner wanted to emphasize that the problem is ambiguous and for him it was enough to show two different procedures leading to two different answers. However, there could be many procedures. Let me just suggest two more procedures that demonstrate the full range of ambiguity.

\begin{enumerate}
\item Suppose that Mr. Smith wants to brag about his sons and will always mention as many as he can. In this case the procedure might be the following:

 If he has two boys, he says, ``I have two boys." If he has one son, he says ``At least one of them is a boy." In this case the answer to the problem is 0.
\item Suppose Mr. Smith doesn't like boys, and wants to de-emphasize the number of boys he has. In this case the procedure might be the following:

If he has two boys, he says, ``At least one of them is a boy." If he has a boy and a girl, he says, ``I am the proud father of a girl." In this case the answer is 1.
\end{enumerate}

I leave it to the reader to invent a procedure for any answer between 0 and 1.

\section{Tuesday-Child Problem}

Now we fast-forward to 2010. Gary Foshee gave a very short talk at the 9th Gathering for Gardner. He said:

\begin{quote}
I have two children. One is a boy born on a Tuesday. What is the probability that I have two boys?
\end{quote}

This is Martin Gardner's Two-Children problem with an extra twist.

Before discussing the solution, let us agree on some basic assumptions:

\begin{enumerate}
\item Sons and daughters are equally probable. This is not exactly true, but it is a reasonable approximation.
\item For our purposes, twins do not exist. Not only is the proportion of twins in the population small, but because they are usually born on the same day, twins might complicate the calculation.
\item All days of the week are equally probable birthdays. While this isn't actually true --- for example, assisted labors are unlikely to be scheduled for weekends --- it is a reasonable approximation.
\end{enumerate}

\subsection{Wrong Solution}

History repeats itself. Just as occurred around the classical Two-Children problem, many mathematicians have been fighting for the wrong solution.  This is their argument:

Each child can be one of two genders and can be born on one of seven days of the week. Thus, gender and the day present 14 equally probable cases for each child. That in turn makes each two-children family belong to one of 196 equally probable cases. When we restrict all possible cases to the given information that one of the children is a boy born on a Tuesday we get 27 equally probable cases. We can divide these cases into several groups, as follows:

\begin{itemize}
\item 7 cases, where the first child is a son born on a Tuesday and the second child is a daughter. 
\item 7 cases, where the second child is a son born on a Tuesday and the first child is a daughter. 
\item 6 cases, where the first child is a son born on a Tuesday and the second child is a son not born on a Tuesday. 
\item 6 cases, where the second child is a son born on a Tuesday and the first child is a son not born on a Tuesday. 
\item 1 case, where both children are sons born on a Tuesday.
\end{itemize}

There is a total of 27 equally probable cases, 13 of which correspond to two sons. Thus the probability must be 13/27.

This wrong solution was widely published. For example, Keith Devlin published it in his column Devlin's Angle in the article ``Probability Can Bite" \cite{Devlin1}. But just like Martin Gardner, he corrected himself in the next article, ``The Problem with Word Problems" \cite{Devlin2}.

\subsection{Procedures}

The ambiguity that Martin Gardner found in the Two-Children problem is likewise present in the Tuesday-Child problem. To resolve ambiguity, we need to specify the procedure by which the information was obtained. Let me discuss four procedures. The calculation is based on 196 equally probable cases for different combinations of gender and the day of the week.

\subsubsection*{Gender-Neutral, Day-of-the-Week-Neutral Procedure}

In this scenario, a father of two children is picked at random. He is instructed to choose a child by flipping a coin. Then he has to provide true information about the chosen child in the following format: ``I have a son/daughter born on Mon/Tue/Wed/Thu/Fri/Sat/Sun." If his statement is, ``I have a son born on a Tuesday," what is the probability that the other child is also a son?

The solution for this procedure is the following: A father has two daughters in 49 cases. Such a father will make the above statement with probability zero. A father has a son and a daughter in 98 cases, and will produce the above statement with a probability of 1/14: with a probability of 1/2 the son is chosen over the daughter and with a probability of 1/7 Tuesday is the birthday. A father has two sons in 49 cases, and he will make the statement with a probability of 1/7. The father of two sons is twice as likely to make the statement as the father of a son and a daughter, but there are half as many such fathers. Thus the probability is 1/2 that the other child is also a son.

This is the most symmetric scenario, which produces the most symmetric answer.

\subsubsection*{Boy-Centered, Day-of-the-Week-Neutral Procedure}

Now let us consider the second scenario. A father of two children is picked at random. If he has two daughters he is sent home and another father picked at random until one is found who has at least one son. If he has one son, he is instructed to provide true information on his son's day of birth. If he has two sons, he has to choose one son at random. His statement will be, ``I have a son born on Mon/Tue/Wed/Thu/Fri/Sat/Sun." If his statement is, ``I have a son born on a Tuesday," what is the probability that the other child is also a son?

The solution is as follows. A father has a son and a daughter in 98 cases, and he will produce the statement above  with a probability of 1/7. If he has two sons (49 cases), the probability of the statement above will likewise be 1/7. The father of two sons is exactly as likely to make the statement as the father of a son and a daughter, but there are half as many such fathers. Thus the probability is 1/3 that the other child is also a son.

This scenario corresponds to the original procedure leading to the first solution of Martin Gardner's Two-Children problem. Unsurprisingly, the answer is the same.

\subsubsection*{Boy-Centered, Tuesday-Centered Procedure}

Now let us consider the third scenario. A father of two children is picked at random. If he doesn't have a son born on a Tuesday, he is sent home and another father is picked at random until one who has a son born on a Tuesday is found. He is instructed to tell you, ``I have a son born on a Tuesday." What is the probability that the other child is also a son?

Here is the solution. A father of a boy and a girl has a son born on a Tuesday in 14 cases. He will make the statement in question with a probability of 1. A father of two sons has a son born on a Tuesday in 13 cases. He too is guaranteed to make the statement. Thus the probability is 13/27 that the other child is a son.

This procedure corresponds to the procedure many mathematicians assume while solving the Tuesday-Child problem. This assumption, made implicitly, is the source of the erroneous solution.

\subsubsection*{Gender-Neutral, Tuesday-Centered Procedure}

For completeness let us consider the fourth scenario. A father of two children is picked at random. If he doesn't have a child who is born on a Tuesday, he is sent home and another father is picked at random until one who has a child born on a Tuesday is found. He is instructed to tell you, ``I have a son/daughter born on a Tuesday." If both of his children were born on Tuesdays, he has to pick one at random. If his statement is, ``I have a son born on a Tuesday," what is the probability that the other child is also a son?

Here is the solution. A father of two daughters will have a child born on a Tuesday in 13 cases. He makes the statement in question with a probability of 0. A father of a boy and a girl has a child  born on a Tuesday in 26 cases. The probability that he makes the statement is 1/2. A father of two sons will have a son born on a Tuesday in 13 cases. The probability that he makes the statement is 1. The father of two sons is twice as likely to make the statement as the father of a son and a daughter, but there are half as many such fathers. Thus the probability is 1/2 that the other child is also a son.

In this procedure we added an additional constraint on the families with two children --- that a child was born on a Tuesday --- that doesn't correlate with gender. Not surprisingly, the answer in this case is 1/2.

\subsection{Correct Solution}

Now let's go back to the original problem. Suppose you meet your friend who you know has two children and he tells you, ``I have a son born on a Tuesday." What is the probability that the other child is also a son?

The problem is under-defined. The solution depends on the reason the father mentions only the son and only the Tuesday.

The funny part of this story is that I, Tanya Khovanova, have two children. And the following statement is true: ``I have a son born on a Tuesday." What is the probability that my other child is a son?

\section{Gender Bias}

What puzzles me is that I've never run into a similar problem about daughters or mothers. I've discussed the original math problem about Mr. Smith with many people many times. But I kept stumbling upon men who passionately defended their wrong solution. When I dug into why their solution was wrong, it appeared that they implicitly assumed that if a man has a daughter and a son, he won't bother talking about his daughter at all.

I've heard this so often that I began to wonder if gender bias wasn't the underlying source of the wrong solution.

Now we have a Tuesday-Child problem. What encourages me is that the most common mistake is that people choose the boy-centered, Tuesday-centered solution. I do not anticipate a strong bias for Tuesday, so perhaps after all it's not a gender bias but simply a mistake.

On second thought, in the father's statement there is a symmetry between genders and days of the week. Why is everyone asking about the gender of the other child, not about the birth day?

I have yet to see people fighting over the following problem.

\begin{quote}
You run into an old friend. He has two children, but you do not know what their gender is. He says, ``I have a son born on a Tuesday." What is the probability that his other child is also born on a Tuesday?
\end{quote}

Real gender equality will be reached when mathematicians start arguing about the probability of the second birthday being Tuesday with the same passion as the probability of the other child being a boy.

\section{Jack}

I have had to defend the solution --- that the problem is ambiguous --- so many times that I invented a fictional opponent, Jack, and here is my imaginary conversation with him.

Jack: The probability that a father with two sons has a son born on a Tuesday is 13/49. The probability that a father with a son and a daughter has a son born on a Tuesday is 1/7. A dad with a son and a daughter is encountered twice as often as a dad with just two sons. Hence, we have probabilities of 13/49 and 14/49, and the probability of the father having a second son is $13/(13+14)$, or 13/27.

Me: What if the problem is about Wednesday?

Jack: It doesn't matter. The particular day in question was random. The answer should be the same: 13/27.

Me: Suppose the father says, ``I have a son born on *day." He mumbles the day, so you do not hear it exactly.

Jack: Well, as the answer is the same for any day, it shouldn't matter. The probability that his other child will also be a son is still 13/27.

Me: Suppose he says, ``I have a son born $\ldots$". So he might have continued and mentioned the day, he might not have. What is the probability?

Jack: We already decided that it doesn't depend on the day, so it shouldn't matter. The probability is still 13/27.

Me: Suppose he says, ``I have a son and I do not remember when he was born." Isn't that the same as just saying, ``I have a son." And by your arguments the probability that his other child is also a son is 13/27.

Jack: Hmm.

Me: Do you remember your calculation? If we denote the number of days in a week as $d$, then the probability that he has another son is $(2d-1)/(4d-1)$. My point is that this probability depends on the number of days in a week. So, if tomorrow we change a week length to another number his probability of having a son changes. Right?

At this point my imaginary conversation stops and I do not know whether I have convinced Jack or not.

\section{Unambiguous Problem}

Now let me give you a variation of the Tuesday-Child problem that is unambiguous and where the answer is 13/27:

\begin{quote}
You pick a random father of two children and ask him, ``Yes or no, do you have a son born on a Tuesday?" Let's make a leap and assume that all fathers know the days of the births of their children and that they answer truthfully. If the answer is yes, what is the probability that the father has two sons?
\end{quote}

The 13/27 argument works perfectly in this case.

The reason this problem stops being ambiguous is that the exact procedure of how you get your information is provided.

\section{Back to Gardner}

Many people I argued with didn't want to listen to me. They referred to Martin Gardner as the final authority to support their wrong solution. Gardner was a great thinker, and he corrected his mistake. I urge those who do not agree with me to trust Martin Gardner and revisit and rethink this  problem together with him.

\section{Acknowledgements}

I am grateful to all my friends and colleagues who discussed the problem with me and supported me in writing about this problem for my blog \cite{blog1, blog2, blog3}. I am especially grateful to Alexey Radul \cite{Radul} and Peter Winkler \cite{Winkler} who contributed their essays on the subject to my blog. I would also like to thank Sue Katz and Julie Sussman, P.P.A., for editing.

%
%
%
%
%


\begin{thebibliography}{9}

\bibitem{Devlin1} Keith Devlin, Probability Can Bite, \url{http://www.maa.org/devlin/devlin_04_10.html}, 2010.

\bibitem{Devlin2} Keith Devlin, The Problem with Word Problems, \url{http://www.maa.org/devlin/devlin_05_10.html}, 2010.

\bibitem{Gardner} Martin Gardner, \emph{The Second Scientific American Book of Mathematical Puzzles and Diversions}, The Chicago University Press. 1987.

\bibitem{blog1} Tanya Khovanova, A Son Born on Tuesday, \url{http://blog.tanyakhovanova.com/?p=221}. 2010.

\bibitem{blog2} Tanya Khovanova, Sons and Tuesdays, \url{http://blog.tanyakhovanova.com/?p=233}. 2010.

\bibitem{blog3} Tanya Khovanova, A Tuesday Quiz, \url{http://blog.tanyakhovanova.com/?p=247}. 2010.

\bibitem{Radul} Alexey Radul, Shannon Entropy Rescues the Tuesday Child, \url{http://blog.tanyakhovanova.com/?p=254}, 2010.

\bibitem{Winkler} Peter Winkler, Conditional Probability and ``He Said, She Said", \url{http://blog.tanyakhovanova.com/?p=234}, 2010.

\end{thebibliography}
\end{document}